\documentclass[12pt]{amsart}
\usepackage{latexsym, amsmath, amscd, amssymb, amsthm} 
\newtheorem{theorem*}{Theorem}[section]

 \newtheorem{lemma*}{Lema}[section] 
  \def\mathbi#1{\textbf{\em #1}}
\begin{document}
\noindent
\title[ Multivariate Poincar\'e series]{  Multivariate Poincar\'e series for  algebras of $SL_2$-invariants}

\author{Leonid Bedratyuk}\address{Khmelnitskiy national university, Instytuts'ka, 11,  Khmelnits'ky, 29016, Ukraine}

\begin{abstract} 
Let $\mathcal{C}_{\mathbi{d}},$ $\mathcal{I}_{\mathbi{d}},$ $\mathbi{d}{=}(d_1,d_2,\ldots, d_n)$ be  the algebras  of join covariants and  joint invariants of the $n$ binary forms of degrees $d_1,d_2,\ldots, d_n.$  Formulas for computation of the multivariate Poincar\'e series $\mathcal{P}(\mathcal{C}_{\mathbi{d}},z_1,z_2,\ldots,z_n,t)$  and  $\mathcal{P}(\mathcal{I}_{\mathbi{d}},z_1,z_2,\ldots,z_n)$  are  found.
 \end{abstract}
\maketitle

{\bf 1. Introduction.} 
Let $V_d$ be the complex vector space  of  binary forms of degree $d$ endowed with the natural action of the special linear group  $G=SL_2.$  Consider the corresponding action of the group  $G$ on  the coordinate rings $\mathbb{C}[V_{\mathbi{d}}]$ and $\mathbb{C}[V_{\mathbi{d}} \oplus \mathbb{C}^2 ],$ where $V_{\mathbi{d}}:=V_{d_1} \oplus V_{d_1} \oplus \cdots \oplus V_{d_n}.$
Denote   by ${\mathcal{I}_{\mathbi{d}}=\mathbb{C}[V_{\mathbi{d}}  ]^{\,G}}$ and by ${\mathcal{C}_{\mathbi{d}}=\mathbb{C}[V_{\mathbi{d}} \oplus \mathbb{C}^2 ]^{\,G}}$ the  subalgebras of   $G$-invariant polynomial functions.
 In the language  of classical invariant theory  the algebras  $\mathcal{I}_{\mathbi{d}}$ and   $\mathcal{C}_{\mathbi{d}}$ are called the algebra  of join invariants and the algebra of join covariants for  the $n$  binary form of  degrees $d_1,d_2, \ldots,d_n$ respectively. The coordinate ring $\mathbb{C}[V_{\mathbi{d}} \oplus \mathbb{C}^2 ]$  can be identified with the algebra of polynomials of the coefficients of the binary forms and of two subsidiary variables, say $X,Y.$  The degree of a covariant    with respect to the variables $X,Y$ is called {\it the order} of the covariant. Every invariant has the order zero. 
 The algebras  $\mathcal{C}_{\mathbi{d}},$ $\mathcal{I}_{\mathbi{d}}$ are  a finitely  generated multigraded algebras under the multidegree-order: 
$$
\mathcal{C}_{\mathbi{d}}=(\mathcal{C}_{\mathbi{d}})_{\mathbi{m},0}+(\mathcal{C}_{\mathbi{d}})_{\mathbi{m},1}+\cdots+(\mathcal{C}_{\mathbi{d}})_{\mathbi{m},j}+ \cdots,
$$
where  each subspace   $(\mathcal{C}_{\mathbi{d}})_{\mathbi{d},j}$ of covariants of multidegree $\mathbi{m}:=(m_1,m_2,\ldots,m_n)$  and order $j$ is   finite-dimensional. The formal power series 
$$
\mathcal{P}(\mathcal{C}_{\mathbi{d}},z_1,z_2,\ldots,z_n,t)=\sum_{\mathbi{m},j=0}^{\infty }\dim((\mathcal{C}_{\mathbi{d}})_{\mathbi{m},j}) z_1^{m_1} z_2^{m_2}\cdots z_n^{m_n} t^j,
$$ 
$$
\mathcal{P}(\mathcal{I}_{\mathbi{d}},z_1,z_2,\ldots,z_n)=\sum_{\mathbi{m}}^{\infty }\dim((\mathcal{I}_{\mathbi{d}})_{\mathbi{m},j}) z_1^{m_1} z_2^{m_2}\cdots z_n^{m_n},
$$

are called the multivariariate Poincar\'e series   of the algebra of   join covariants  $\mathcal{C}_{\mathbi{d}}.$
It is clear that the series  $\mathcal{P}(\mathcal{C}_{\mathbi{d}},z,z,\ldots,z,0)$ is the Poincare series of the algebra $\mathcal{I}_{\mathbi{d}}$  and the series $\mathcal{P}(\mathcal{C}_{\mathbi{d}},z,z,\ldots,z,1)$ is the Poincare series of the algebra $\mathcal{C}_{\mathbi{d}}$  with respect to the usual grading of the algebra under  degree.
 The algebra of covariants $\mathcal{C}_d$ is Cohen-Macaulay. It  implies  that its multivariate Poincar\'e series is the  power series expansion  of a  rational function of the  variables $z_1,\ldots,z_n,t.$  We consider here the problem of
computing efficiently this rational function.

 Calculating  the  Poincar\'e series   of the  algebras of invariants and covariants was an  important  object of research  in  invariant theory in  the 19th century.
For  the case one binary form   $d\leq 10,$ $d=12$ the  bivariate  series $\mathcal{P}(\mathcal{C}_d,z,t)$  were calculated by Sylvester  and Franklin, see  in \cite{SF, Sylv-12} the big tables  named them as Generating Functions for covariants, reduced form. All those calculations are correct  up to $d=6.$ In  \cite{Br}-\cite{BC2} the Poincar\'e series of algebras of   joint  covariants of two and three binary form of small degrees  are  calculated.

It is well known that  there is an isomorphism between $\mathcal{C}_{\mathbi{d}}$  and $\ker \mathcal{D}_{\mathbi{d}},$ where $\mathcal{D}_{\mathbi{d}}$  is the Weitzenb\"ock derivarion, see  details in \cite{CIT-LND_I}.  Therefore, we   have $\mathcal{P}(\mathcal{C}_{\mathbi{d}},z_1,z_2,\ldots,z_n,t)=\mathcal{P}(\ker \mathcal{D}_{\mathbi{d}},z_1,z_2,\ldots,z_n,t)$.

 In the paper we  offer  a  Cayley-Sylvester type formula  for calculating of $\dim (\mathcal{C}_{\mathbi{d}})_{\mathbi{m},i}$   and a  formulas  for calculation of $\mathcal{P}(\mathcal{C}_{\mathbi{d}},z_1,z_2,\ldots,z_n,t)$ and $\mathcal{P}(\mathcal{I}_{\mathbi{d}},z_1,z_2,\ldots,z_n,t).$ By using the  formula  and Xin's package for MacMahon Partitions analisys, the multivariate Poincar\'e series $\mathcal{P}(\mathcal{C}_{\mathbi{d}},z_1,z_2,\ldots,z_n,t),$ $\mathcal{P}(\mathcal{I}_{\mathbi{d}},z_1,z_2,\ldots,z_n)$ for small $\mathbi{d}$ is calculated. Also, we offer a Maple package {\tt Poincare\_series.}


{\bf 2. Cayley-Sylvester type formula}


 To begin, we give a proof of  Cayley-Sylvester type  formula for the dimension of multi-graded subspaces   of   the algebra $\mathcal{C}_{\mathbi{d}}.$

Let  $V_{d_k}^{(k)}=\langle v_0^{(k)},v_1^{(k)},...,v_{d_k}^{(k)} \rangle,$ $\dim V_d^{(k)}=d_k+1,$ $k=1,\ldots,n$ be the set of $n$ standard  irreducible  representations of the complex  Lie algebra  $\mathfrak{sl_{2}}.$
The basis elements    $ \left( \begin{array}{ll}  0\, 1 \\ 0\,0 \end{array} \right),$ $ \left( \begin{array}{ll}  0\, 0 \\ 1\,0 \end{array} \right)$, $ \left( \begin{array}{ll}  1 &  \phantom{-}0 \\  0 &-1 \end{array} \right)$ of the algebra    $\mathfrak{sl_{2}}$ act on    $V_{d_k}$  by the derivations  $D_1, D_2, E:$ 
$$
D_1 \left(v_i^{(k)}\right)=i\, v_{i-1}^{(k)},  D_2\left(v_i^{(k)}\right)=(d-i)\,v_{i+1}^{(k)}, E\left(v_i^{(k)}\right)=(d-2\,i)\,v_i^{(k)}.
$$
The action of   $\mathfrak{sl_{2}}$  is extended to an action on the symmetrical algebra  $S(V_{\mathbi{d}})$ in the natural way. 

Let   $\mathfrak{u}_{2}$ be  the maximal unipotent subalgebra of $\mathfrak{sl}_{2}.$ The algebra  $\mathcal{S}_{\mathbi{d}},$ defined by 
$$
\mathcal{S}_{\mathbi{d}}:= \displaystyle{ S(V_{\mathbi{d}})^{\mathfrak{u_{2}}}}=\{ v \in S(V_{\mathbi{d}})|  D_1(v)=0 \},
$$
is called the {\it  algebra of joint semi-invariants}  of the $n$ binary form  of degrees $d_1,d_2,\ldots,d_n.$ For any element $v \in \mathcal{S}_{\mathbi{d}}$ a natural number $s$ is called {\it the  order} of the semi-invariant $v$ if the number $s$ is the smallest natural number such that \begin{equation*}D_2^s(v) \ne 0, D_2^{s+1}(v) = 0.\end{equation*}
It is  clear that any semi-invariant  $v \in \mathcal{S}_{\mathbi{d}}$  of  order $i$ is the highest weight vector  for an  irreducible $\mathfrak{sl_{2}}$-module   of the dimension $i+1$ in $S(V_{\mathbi{d}}).$

In  \cite{CIT-LND_I} we proved  that  the algebras of joint covariants and the algebra of joint semi-invariants  are  isomorphic. Furthermore, the order of elements  is preserved through the isomorphism.  Thus, it is  enough to compute the Poincar\'e  series of the algebra of semi-invariants  $\mathcal{S}_{\mathbi{d}}.$ 

The algebra  $S(V_{\mathbi{d}})$  is $\mathbb{N}^{n}$-graded 
$$
S(V_{\mathbi{d}})=S^{(0,0,\ldots,0)}(V_{\mathbi{d}})+S^{(1,0,\ldots,0)}(V_{\mathbi{d}})+\cdots +S^{\mathbi{m}}(V_{\mathbi{d}})+\cdots,
$$
and each   $S^{\mathbi{m}}(V_{\mathbi{d}}),$ $\mathbi{m} \in \mathbb{N}^{n} $ is a completely  reducible
 representation of the Lie algebra  $\mathfrak{sl_{2}}.$

Let  $V_k$ be  the standard irreducible    $\mathfrak{sl_{2}}$-module, $\dim V_k=k+1.$ Then, the following primary decomposition  holds
$$
S(V_{\mathbi{d}})^{\mathbi{m}} \cong \gamma_\mathbi{d}({\mathbi{m}};0) V_0+\gamma_\mathbi{d}({\mathbi{m}};1) V_1+ \cdots +\gamma_\mathbi{d}({\mathbi{m}};md^*) V_{md^*},  \eqno{(1)}
$$
here  $md^*:=\max(m_1\,d_1,m_2\,d_2,\ldots m_n\,d_n)$ and $\gamma_\mathbi{d}({\mathbi{m}};i)$ is  the  multiplicity of the representation  $V_k$  in the decomposition of  $S(V_{\mathbi{d}})^{\mathbi{m}}.$ On the other hand, the multiplicity  $\gamma_\mathbi{d}({\mathbi{m}};i)$  is  equal to the number of linearly independent homogeneous joint semi-invariants of  the multidegree $\mathbi{m}$   and the order $i.$
In particular, the number of linearly  independent joint invariants of the multidegree $\mathbi{m}$ is  equal to  $\gamma_\mathbi{d}({\mathbi{m}};0).$  These arguments prove

{\bf Lemma 1.}
$$
 \dim (S_{\mathbi{d}})_{\mathbi{m},i}=\gamma_{{\mathbi{d}}}({\mathbi{m}};i).
$$

Denote by  $\Lambda_{W}$ the set of all weights  of a representation  $W,$  for instance  $$\Lambda_{V_d}=\{-d, -d+2, \ldots,d-2, d \}.$$
 
Consider the    set of variables:  $v^{(1)}_{0},v^{(1)}_{1},\ldots, v^{(1)}_{d_1},$ $v^{(2)}_{0},v^{(2)}_{1},\ldots, v^{(2)}_{d_2},$ $\ldots,$ $v^{(n)}_{0},v^{(n)}_{1},\ldots, v^{(n)}_{d_n}.$
The character  $ {\rm Char}\left(S(V_{\mathbi{d}})^{\mathbi{m}}\right)$ of the representation  $S(V_{\mathbi{d}})^{\mathbi{m}},$ see  \cite{FH}, equals   $$H_\mathbi{m}(q^{-d_1},q^{-d_1+2},\ldots,q^{d_1},q^{-d_2},q^{-d_2+2},\ldots,q^{d_2},\ldots,q^{-d_n},q^{-d_n+2},\ldots,q^{d_n}),$$   where  $H_{\mathbi{m}}(v^{(1)}_{0},v^{(1)}_{1},\ldots, v^{(1)}_{d_1},\ldots,v^{(n)}_{0},v^{(n)}_{1},\ldots, v^{(n)}_{d_n})$ is  the complete symmetrical function     
$$
\begin{array}{l}
\displaystyle H_\mathbi{m}(v^{(1)}_{0},v^{(1)}_{1},\ldots, v^{(1)}_{d_1},\ldots,v^{(n)}_{0},v^{(n)}_{1},\ldots, v^{(n)}_{d_n})=\\
\\
\displaystyle =\sum_{|\alpha^{(1)}|=m_1,\ldots,|\alpha^{(n)}|=m_n} (v^{(1)}_{0})^{\alpha_0^{(1)}}(v^{(1)}_{1})^{\alpha_1^{(1)}}\ldots (v^{(1)}_{d_1})^{\alpha_{d_1}^{(1)}}\cdots (v^{(n)}_{0})^{\alpha_0^{(n)}} \ldots (v^{(1)}_{d_n})^{\alpha_{d_n}^{(n)}},
\end{array}
$$
where  $\displaystyle |\alpha^{(k)}|:=\sum_{i=0}^{d_i}\alpha_i^{(k)}.$

Replacing   $v_i^{(k)}$ by $q^{d_k-2\,i},$   we  obtain the specialized expression for the character   $S(V_{\mathbi{d}})^{\mathbi{m}},$ namely 
$$
\begin{array}{c}
\displaystyle {\rm Char}(S(V_{\mathbi{d}})^{\mathbi{m}})= \\
\\
\displaystyle = \sum_{|\alpha^{(1)}|=m_1,\ldots,|\alpha^{(n)}|=m_n} q^{d_1|\alpha^{(1)}|+\ldots+d_n|\alpha^{(n)}| -2\left(\alpha_1^{(1)}+2\alpha_2^{(1)}+\cdots + d_1\, \alpha_{d_1}^{(1)}\right)-\ldots-2\left(\alpha_1^{(n)}+2\alpha_2^{(n)}+\cdots + d_n\, \alpha_{d_n}^{(n)}\right)}=\\
\\
\displaystyle =\sum_{i=-md^*}^{md^*} \omega_\mathbi{d}({\mathbi{m}};i) q^{i},
\end{array}
$$
here    $\omega_\mathbi{d}({\mathbi{m}};i)$  is the number of non-negative integer solutions of the following  system of equations:
$$
\left \{
\begin{array}{l}
d_1|\alpha^{(1)}|+\ldots+d_n|\alpha^{(n)}| -2\left(\alpha_1^{(1)}+2\alpha_2^{(1)}+\cdots + d_1\, \alpha_{d_1}^{(1)}\right)-\\
-\ldots-2\left(\alpha_1^{(n)}+2\alpha_2^{(n)}+\cdots + d_n\, \alpha_{d_n}^{(n)}\right)=i, \\
\\
|\alpha^{(1)}|=m_1,\\
|\alpha^{(2)}|=m_2,\\
\ldots \\
|\alpha^{(s)}|=m_n.\\
\end{array}
\right. 
 \eqno{(2)}
$$

We can summarize what we have shown so far in  

{\bf Theorem 1.} 
$$
 \dim (\mathcal{S}_{{\mathbi{d}}})_{\mathbi{m},i}=\omega_{\mathbi{d}}({\mathbi{m}};i)-\omega_{\mathbi{d}}({\mathbi{m}};i+2).
$$

Let us transform  the system  $(2)$  to the  form:
$$ 
\left \{
\begin{array}{l}
d_1\alpha_0^{(1)}+(d_1-2)\alpha_1^{(1)}+(d_1-4)\alpha_2^{(1)}+\cdots + (-d_1)\, \alpha_{d_1}^{(1)}+\cdots +\\ \\+d_s\alpha_0^{(s)}+(d_s-2)\alpha_1^{(s)}+(d_s-4)\alpha_2^{(s)}+\cdots + (-d_s)\, \alpha_{d_s}^{(s)}=i, \\
\\
|\alpha^{(1)}|=m_1,\\
|\alpha^{(2)}|=m_2,\\
\ldots \\
|\alpha^{(s)}|=m_s.\\
\end{array}
\right. 
$$
 It well-known  that  the number
 $\omega_{\mathbi{d}}({\mathbi{m}};i)$ of   non-negative integer solutions of the above  system
is equal to the coefficient of  $\displaystyle z_1^{m_1} z_2^{m_2}\cdots z_n^{m_n} t^i $ of the  generating function
$$
\begin{array}{l}
f_{{\mathbi{d}}}(z_1,z_2,\ldots,z_n,t)=\displaystyle \frac{1}{\displaystyle \prod_{k=1}^n \prod_{j=0}^{d_k} (1-z_k t^{d_k-2\,j})}.\\
\end{array}
$$
Denote it in such a way:  $\omega_{\mathbi{d}}({\mathbi{m}};i):=\left[  z^{{\mathbi{m}}} t^i\right]f_{{\mathbi{d}}}(z_1,z_2,\ldots,z_n,t).$ 
Observe that $$ f_{{\mathbi{d}}}(z_1,z_2,\ldots,z_n,t)=f_{{\mathbi{d}}}(z_1,z_2,\ldots,z_n,t^{-1}).$$

The  following statement holds

{\bf Theorem 2.}
$$
\begin{array}{l}
 \dim (\mathcal{S}_{{\mathbi{d}}})_{\mathbi{m},i}=[z^{\mathbi{m} }t^{d_1 m_1+d_2 m_2+\cdots + d_n m_n-\,i} ](1-t^{2})f_{{\mathbi{d}}}(z_1 t^{d_1},z_2 t^{d_2},\ldots,z_n t^{d_n},t).
\end{array}
$$

We  replaced  $z_k$  with   $z_k t^{d_k},k=1,..,n$  to avoid  of a negative powers  of $t$  in the denominator of the function $f_{{\mathbi{d}}}(\mathbi{z},t).$

{\bf 3.Computing the Poincar\'e series by MacMahon partition analysis}

 Let us present    formulas for the Poincar\'e  series  for the algebras $\mathcal{S}_{\mathbi{d}} \cong \mathcal{C}_{\mathbi{d}}$  and $\mathcal{I}_{\mathbi{d}}.$  
Consider the $\mathbb{C}$-algebra $\mathbb{C}[[z_1,z_2,\ldots, z_n,t]]$   of  the formal   power series.
Define  $\mathbb{C}$-linear function
$$ \Psi_{\mathbi{d}}:\mathbb{C}[[z_1,z_2,\ldots, z_n,t]] \to \mathbb{C}[[z_1,z_2,\ldots, z_n,t]],$$  in the  following  way:
$$
\Psi_{\mathbi{d}}\left(\sum_{{\mathbi{i}},j=0}^{\infty} a_{{\mathbi{i}},j}\, z^{{\mathbi{i}}} t^j\right)=\sum_{d_1 i_1+\cdots d_n i_n-j\geq 0} a_{{\bf{i}},j} z^{{\mathbi{i}}} t^{d_1 i_1+\cdots d_n i_n-j},
$$
here $a_{{\mathbi{i}},j}:=a_{i_1,i_2,\ldots,i_n,j} \in \mathbb{C}.$

The main idea of the calculations of the paper  is that 
the  multivariate Poincar\'e series   can be expressed  in terms of the function $ \Psi_{\mathbi{d}}.$ The following simple but important statement  holds:

{\bf Lemma 2.}    
$
\mathcal{P}(\mathcal{S}_{\mathbi{d}},z_1,z_2,\ldots,z_n,t)=\Psi_{\mathbi{d}} \left((1-t^{2})f_{{\mathbi{d}}}(z_1 t^{d_1},z_2 t^{d_2},\ldots,z_n t^{d_n},t)\right).
$

To compute  the multivariate Poincar\'e  series $\mathcal{P}(\mathcal{C}_{\mathbi{d}},z_1,\ldots z_n,t)$ and $\mathcal{P}(\mathcal{I}_{\mathbi{d}},z_1,\ldots z_n)$  we use the well-known   MacMahon's Omega operators $\underset{\scriptscriptstyle \geq 0}{\Omega},$ $\underset{\scriptscriptstyle = 0}{\Omega}$ which act on the  Laurent series 
$$
 \sum_{k_1=-\infty}^{\infty} \cdots  \sum_{k_s=-\infty}^{\infty}\sum_{\alpha=-\infty}^{\infty} a_{k_1,k_2,\ldots,k_s,\alpha} z_1^{k_1}z_2^{k_2}\cdots z_s^{k_s} \lambda^{\alpha}
$$
by 
$$
\underset{\scriptscriptstyle \geq 0}{\Omega}\sum_{k_1=-\infty}^{\infty} \cdots  \sum_{k_s=-\infty}^{\infty}\sum_{\alpha=-\infty}^{\infty} a_{k_1,\ldots,k_s,\alpha} z_1^{k_1}\cdots z_s^{k_s} \lambda^{\alpha}=\sum_{k_1=0}^{\infty} \cdots  \sum_{k_s=0}^{\infty} \sum_{\alpha=0}^{\infty} a_{k_1,\ldots,k_s,\alpha} z_1^{k_1}\cdots z_s^{k_s} \lambda^{\alpha},
$$
and
$$
\underset{\scriptscriptstyle = 0}{\Omega}\sum_{k_1=-\infty}^{\infty} \cdots  \sum_{k_s=-\infty}^{\infty}\sum_{\alpha=-\infty}^{\infty} a_{k_1,\ldots,k_s,\alpha} z_1^{k_1}\cdots z_s^{k_s} \lambda^{\alpha}=\sum_{k_1=0}^{\infty} \cdots  \sum_{k_s=0}^{\infty} a_{k_1,\ldots,k_s,\alpha} z_1^{k_1}\cdots z_s^{k_s}
$$

The following statement holds:

{\bf Theorem  2.}
$$
\mathcal{P}(\mathcal{C}_{\mathbi{d}},z_1,\ldots z_n,t)=\underset{\scriptscriptstyle \geq 0}{\Omega}f_{\mathbi{d}}\left(z_1 (t\lambda)^{d_1},z_2 (t\lambda)^{d_2},\ldots,z_s (t\lambda)^{d_s},\dfrac{1}{t\lambda}\right),
$$
$$
\mathcal{P}(\mathcal{I}_{\mathbi{d}},z_1,\ldots z_n)=\underset{\scriptscriptstyle = 0}{\Omega}f_{\mathbi{d}}\left(z_1 (t\lambda)^{d_1},z_2 (t\lambda)^{d_2},\ldots,z_s (t\lambda)^{d_s},\dfrac{1}{t\lambda}\right).
$$

To calculate the multigraded Poincar\'e series  we have implemented the above results in  our Maple package  {\tt Poincare\_series}. The last version of the package  can be dowloaded from the web site http://sites.google.com/site/bedratyuklp.

The  package's  procedure {\tt  MULTIVAR\_COVARIANTS}~$([d_1,d_2,\ldots,d_n])$ calculates  the multivariate Poincar\'e  series 
 $\mathcal{P}(\mathcal{C}_{\mathbi{d}},z_1,z_2,\ldots,z_n,t).$ 
Below are some examples:
\begin{gather*}
\mathcal{P}(\mathcal{C}_{(1,1)},z_1,z_2,t)={\frac {1}{ \left( 1-z_{{2}}t \right)  \left( 1-z_{{1}}t \right) 
 \left( 1-z_{{1}}z_{{2}} \right) }},\\ \mathcal{P}(\mathcal{C}_{(1,2)},z_1,z_2,t)={\frac {1+z_{{1}}z_{{2}}t}{ \left( 1-z_{{2}}{t}^{2} \right)  \left( 
1-{z_{{2}}}^{2} \right)  \left( 1-z_{{1}}t \right)  \left( 1-{z_{{1}
}}^{2}z_{{2}} \right) }},\\
\mathcal{P}(\mathcal{C}_{(2,2)},z_1,z_2,t)={\frac {1+z_{{1}}z_{{2}}{t}^{2}}{ \left( 1-z_{{1}} \right)  \left( 
1-z_{{2}}{t}^{2} \right)  \left(1- {z_{{2}}}^{2} \right)  \left( 1-z_
{{1}}{t}^{2} \right)  \left( 1-z_{{2}}z_{{1}} \right)  \left(1- z_{{1}}
 \right) }}.
\end{gather*}

Also, the procedure {\tt  MULTIVAR\_INVARIANTS}~$([d_1,d_2,\ldots,d_n]),$ calculates the  multivariate series  $\mathcal{P}(\mathcal{I}_{\mathbi{d}},z_1,z_2,\ldots,z_n),$ for instance
\begin{gather*}
\mathcal{P}(\mathcal{I}_{(1,1)},z_1,z_2)=\frac{1}{ 1-z_{{1}}z_{{2}} },
\mathcal{P}(\mathcal{I}_{(1,3)},z_1,z_2)={\frac {1+{z_{{2}}}^{2}{z_{{1}}}^{2}-z_{{2}}z_{{1}}}{ \left(1-{z_{{2}}}^{4} \right)  \left( 1-{z_{{1}}}^{3}z_{{2}} \right)  \left(1-z_
{{2}}z_{{1}} \right) }},\\
\end{gather*}

\end{document}